\documentclass{amsart}
\usepackage{amssymb,latexsym}
\usepackage[latin1]{inputenc}
\usepackage[T1]{fontenc}
\usepackage{epsf,graphicx}
\usepackage{psfrag}

\theoremstyle{plain}
\newtheorem{lemma}{Lemma}
\newtheorem{proposition}[lemma]{Proposition}
\newtheorem{example}[lemma]{Example}
\newtheorem{remark}[lemma]{Remark}
\newtheorem{theorem}[lemma]{Theorem}

\theoremstyle{definition}
\newtheorem{definition}[lemma]{Definition}

\theoremstyle{remark}
\newtheorem{notation}[lemma]{Notation}

\begin{document}
\title{The Stokes Phenomenon in the Confluence of the Hypergeometric Equation Using Riccati Equation}
\author{Caroline Lambert and Christiane Rousseau}
\address{Département de mathématiques et de statistique\\Université de Montréal\\C.P. 6128, Succursale Centre-ville, Montréal (Qc), H3C 3J7, Canada}
\email{lambert@dms.umontreal.ca, rousseac@dms.umontreal.ca}
\thanks{Research supported by NSERC and FQRNT in Canada}
\keywords{Hypergeometric equation, confluence, Stokes phenomenon, divergent series, analytic continuation, summability, monodromy, confluent hypergeometric equation, Riccati equation} 
\date{\today}
\begin{abstract}
In this paper we study the confluence of two regular singular points of the hypergeometric equation into an irregular one. We study the consequence of the divergence of solutions at the irregular singular point for the unfolded system. Our study covers a full neighborhood of the origin in the confluence parameter space. In particular, we show how the divergence of solutions at the irregular singular point explains the presence of logarithmic terms in the solutions at a regular singular point of the unfolded system. For this study, we consider values of the confluence parameter taken in two sectors covering the complex plane. In each sector, we study the monodromy of a first integral of a Riccati system related to the hypergeometric equation. Then, on each sector, we include the presence of logarithmic terms into a continuous phenomenon and view a Stokes multiplier related to a 1-summable solution as the limit of an obstruction that prevents a pair of eigenvectors of the monodromy operators, one at each singular point, to coincide.
\end{abstract}

\maketitle

\pagestyle{myheadings}\markboth{C.Lambert, C.Rousseau}{Confluence of the hypergeometric equation}

\section{Introduction}

The hypergeometric differential equation arises in many problems of mathematics and physics and is related to special functions. It is written
\begin{equation}\label{E:hypergeo}
X (1-X)\, v''(X) + \{ c-(a+b+1) X \} \, v'(X)-a b \, v(X) = 0.
\end{equation}

More precisely, any linear equation of order two ($y''(z)+p(z)y'(z)+q(z)y(z)=0$) with three regular singular points can be transformed into the hypergeometric equation by a change of variables of the form $y=f(z)v$ and a new independant variable $X$ obtained from $z$ by a Möbius transformation (see for example \cite{eR64}). 

The confluent hypergeometric equation with a regular singular point at $z=0$ and an irregular one at $z=\infty$  is often written in the form 
\begin{equation} \label{eq hyp confluente originale}
z u''(z)+(c'-z)u'(z)- a' u(z)=0.
\end{equation}

Solutions of this equation at the irregular point $z=\infty$ are in general divergent and always 1-summable. C. Zhang (\cite{cZ94} and \cite {cZ96}) and J.-P. Ramis \cite{jR89} showed that the Stokes multipliers related to the confluent equation can be obtained from the limits of the monodromy of the solutions of the nonconfluent equation (\ref{E:hypergeo}). They assumed that the bases of solutions of (\ref{E:hypergeo}) around the merging singular points ($z=b$ and $z=\infty$) never contain logarithmic terms and they described the phenomenon using two types of limits: first with $\Im(b) \to \infty$, then with $\Re(b) \to \infty$ on the subset $b=b_0 + \mathbb{N}$ for $b_0 \in \mathbb{C}$. They also proved the uniform convergence of the solutions on all compact sets in the case $\Im b \to \infty$.

In this paper, we propose a different approach : we describe the phenomenon in a whole neighborhood of values of the confluence parameter, but we are forced to cover the neighborhood with two sectors on which the presentations are different. We are then able to explain the presence of the logarithmic terms: they occur precisely for discrete values of the confluence parameter when we unfold a confluent equation with at least one divergent solution. On each sector, each divergent solution explains the presence of logarithmic terms at one of the unfolded singular points. The occurrence of logarithmic terms, a discrete phenomenon, is embedded into a continuous phenomenon valid on the whole sector. 

To help understanding the phenomenon, we give a translation of the hypergeometric equation in terms of a Riccati system in which two saddle-nodes are unfolded with a parameter $\epsilon$. The parameter space is again covered with two sectors $S^±$. For this Riccati system, we consider on each sector $S^±$ of the parameter space a first integral which has a limit when $\epsilon \to 0$, written in the form $I^{\epsilon^±}(x,y)=H^{\epsilon^±}(x) \frac{y-\rho_1(x,\epsilon)}{y-\rho_2(x,\epsilon)}$ where $y=\rho_1(x,\epsilon)$ and $y=\rho_2(x,\epsilon)$ are analytic invariant manifolds of singular points and, for $\epsilon =0$, center manifolds of the saddle-nodes. Then, when we calculate the monodromy of one of these first integrals, we can separate it into two parts: a continuous one which has a limit when $\epsilon \to 0$ inside the sector $S^±$ and a wild one which has no limit but which is linear. The wild part is independent of the divergence of the solutions and present in all cases. The divergence of $\rho_1(x,0)$ corresponds to the analytic invariant manifold of one singular point being ramified at the other in the unfolding of one saddle-node. For particular values of $\epsilon$ for which one singular point is a resonant node, this forces the node to be nonlinearisable (i.e. to have a nonzero resonant monomial), in which case logarithmic terms appear in $I^{\epsilon^±}$. This is called the parametric resurgence phenomenon in \cite{cR05}. The divergence of $\rho_2(x,0)$ corresponds to a similar phenomenon with the pair of singular points coming from the unfolding of the other saddle-node. Finally, we translate our results in the case of a universal deformation.

\section{Solutions of the hypergeometric equations}\label{S:div ramifi}

In this paper, we study the confluence of the singular points $0$ and $1$; the confluent hypergeometric equation has an irregular singular point at the origin. We make the change of variables $X=\frac{x}{\epsilon}$ in (\ref{E:hypergeo}) to bring the singular point at $X=1$ to a singular point at $x=\epsilon \ne 0$. We consider small values of $\epsilon$ and we limit the values of $c$ to 
\begin{equation}\label{E:c}
c=1-\frac{1}{\epsilon}.
\end{equation}
Let $v(\frac{x}{\epsilon})$ be denoted  by $w(x)$. Then (\ref{E:hypergeo}) becomes
\begin{equation}\label{E:hypergeo e}
x(x-\epsilon) \, w''(x) + \{ 1-\epsilon +(a+b+1)x \} \, w'(x)+a b \, w(x) = 0.
\end{equation}
We will then let $\epsilon \to 0$. We want to study what happens in a neighborhood of $\epsilon =0$. The confluence parameter $\epsilon$ will be taken in two sectors, the union of which is a small pointed neighborhood of the origin in the complex plane. 

\begin{remark}
Although not explicitly written, our study is still valid if we let $a(\epsilon)$ and $b(\epsilon)$ be analytic functions of $\epsilon$.
\end{remark}

\begin{definition}
Given $\gamma \in (0,\frac{\pi}{2})$ fixed, we define
\begin{itemize} 
\item $S^+=\{ \epsilon \in \mathbb{C} \, : \, 0<|\epsilon|<r(\gamma), \, \arg(\epsilon) \in (-\pi+\gamma,\pi -\gamma ) \}$,
\item $S^-=\{ \epsilon \in \mathbb{C} \, : \, 0<|\epsilon|<r(\gamma), \, \arg(\epsilon) \in (\gamma ,2 \pi -\gamma) \}$.
\end{itemize}

\begin{remark}$\gamma$ can be chosen arbitrary small, but $r(\gamma)$ will depend on $\gamma$ and $r(\gamma) \to 0$ as $\gamma \to 0$. In particular, we will ask $a+b+\frac{1}{\epsilon} \notin -\mathbb{N}$, $a+\frac{1}{\epsilon} \notin -\mathbb{N}$ and $b+\frac{1}{\epsilon} \notin -\mathbb{N}$ on $S^+$ and $2-a-b-\frac{1}{\epsilon} \notin -\mathbb{N}$, $a-\frac{1}{\epsilon} \notin -\mathbb{N}$ and $b-\frac{1}{\epsilon} \notin -\mathbb{N}$ on $S^-$ (in this paper $\mathbb{N}=\{0,1,...\}$).
\end{remark} 
\end{definition}

\subsection{Bases for the solutions of the hypergeometric equation (\ref{E:hypergeo e}) at the regular singular points $x=0$ and $x=\epsilon$} \label{S:sol hyp}
The fundamental group of $\mathbb{C} \backslash \{0,\epsilon \}$ based at an ordinary point acts on a solution (valid at this base point) by giving its analytic continuation at the end of a loop. In this way we have monodromy operators around each singular point. We can extend it to act on any function of solutions.

\begin{notation}
The monodromy operator $M_0$ (resp. $M_\epsilon$) is the one associated to the loop which makes one turn around the singular point $x=0$ (resp. $x=\epsilon$) in the positive direction (and which does not surround any other singular point). In this paper, since we use bases of solutions whose Taylor series are convergent in a disk of radius $\epsilon$ centered at a singular point, it will be useful to define $M_0$ (resp. $M_\epsilon$) with the fundamental group based at a point belonging to the line joining $-\epsilon$ and $0$ (resp. $\epsilon$ and $2\epsilon$).
\end{notation}

As the hypergeometric equation is linear of second order, the space of solutions is of dimension $2$. Given a basis for the space of solutions, the monodromy operator $M_0$ (resp. $M_\epsilon$) acting on this basis is linear and is represented by a two-dimensional matrix.

As elements of a basis $\mathcal{B}_0$ (resp. $\mathcal{B}_\epsilon$) around the singular point $x=0$ (resp. $x=\epsilon$), it is classical to use solutions which are eigenvectors of the monodromy operator $M_0$ (resp. $M_\epsilon$) whenever these solutions exist. However, none of these bases is defined on the whole of a sector $S^+$ or $S^-$. This is why we later switch to mixed bases. C. Zhang (\cite{cZ94} and \cite{cZ96}) also used mixed bases but he has not pushed the study as far as we do. 

\begin{definition}
The hypergeometric series $\, _kF_j(a_1,a_2,...a_k,c_1,c_2,...,c_j;x)$ is defined by 
\begin{equation}
\, _kF_j(a_1,a_2,...a_k,c_1,c_2,...,c_j;x)=1+\sum_{n=1}^{\infty} \frac{(a_1)_{n}(a_2)_{n}...(a_k)_{n}}{(c_1)_{n}(c_2)_{n}...(c_j)_{n}n!}x^{n} 
\end{equation}
with
\begin{equation}
\begin{cases}
(a)_0=1 \\
(a)_n=a(a+1)(a+2)...(a+n-1)
\end{cases}
\end{equation}
and for $c_1,...,c_j \notin -\mathbb{N}$.
\end{definition}

A basis $\mathcal{B}_0=\{w_1(x), w_2(x) \}$ of solutions of (\ref{E:hypergeo e}) around the singular point $x=0$ is well known (see \cite{yL69} for details):
\begin{align}\label{E:w1(x) w2(x)}
\begin{cases}
w_1(x)&= \, _2F_1(a,b,1-\frac{1}{\epsilon};\frac{x}{\epsilon}) \\
&=(1-\frac{x}{\epsilon})^{1-\frac{1}{\epsilon}-a-b}\, _2F_1(1-\frac{1}{\epsilon}-a,1-\frac{1}{\epsilon}-b,1-\frac{1}{\epsilon};\frac{x}{\epsilon}), \\
w_2(x)&=(\frac{x}{\epsilon})^{\frac{1}{\epsilon}} \, _2F_1(a+\frac{1}{\epsilon},b+\frac{1}{\epsilon},1+\frac{1}{\epsilon};\frac{x}{\epsilon})\\
&=(\frac{x}{\epsilon})^{\frac{1}{\epsilon}}(1-\frac{x}{\epsilon})^{1-\frac{1}{\epsilon}-a-b}\, _2F_1(1-a,1-b,1+\frac{1}{\epsilon};\frac{x}{\epsilon}).
\end{cases}
\end{align}
The solution $w_1(x)$ exists if $1-\frac{1}{\epsilon} \notin -\mathbb{N}$ whereas $w_2(x)$ exists if $1+\frac{1}{\epsilon} \notin -\mathbb{N}$.

Similarly, a basis $\mathcal{B}_\epsilon= \{w_3(x), w_4(x) \}$ of solutions of (\ref{E:hypergeo e}) around the singular point $x=\epsilon$ is given by:
\begin{align}\label{E:w3(x) w4(x)}
\begin{cases}
w_3(x)&= \, _2F_1(a,b,a+b+\frac{1}{\epsilon};1-\frac{x}{\epsilon}), \\
w_4(x)&=(\frac{x}{\epsilon})^{\frac{1}{\epsilon}}(1-\frac{x}{\epsilon})^{1-\frac{1}{\epsilon}-a-b} \, _2F_1(1-a,1-b,2-\frac{1}{\epsilon}-a-b;1-\frac{x}{\epsilon}).
\end{cases}
\end{align}
The solution $w_3(x)$ exists if $a+b+\frac{1}{\epsilon} \notin -\mathbb{N}$ whereas $w_4(x)$ exists if $2-\frac{1}{\epsilon}-a-b \notin -\mathbb{N}$. 

In particular, $w_2(x)$ and $w_3(x)$ exist for all $\epsilon \in S^+$ and $w_1(x)$ and $w_4(x)$ exist for all $\epsilon \in S^-$, provided $r(\gamma)$ is sufficiently small.

Traditionally, in order to get a basis when $1-\frac{1}{\epsilon} \in -\mathbb{N}$, $a \notin -\mathbb{N}$ and $b \notin -\mathbb{N}$ (resp.  $2-\frac{1}{\epsilon}-a-b \in -\mathbb{N}$, $1-a \notin -\mathbb{N}$ and $1-b \notin -\mathbb{N}$), the solution $w_1(x)$ in $\mathcal{B}_0$ (resp. $w_4(x)$ in $\mathcal{B}_\epsilon$) is replaced by some other solution $\tilde w_1(x)$ (resp. $\tilde w_4(x)$) which contains logarithmic terms. The converse is true if $\epsilon \in S^+$ is sufficiently small. Similarly, we have $\tilde w_2(x)$ and $\tilde w_3(x)$ for specific value of $\epsilon$ in $S^-$ (see for example \cite{eG36}).

The problem with this approach is that the basis $\mathcal{B}_0=\{w_1(x), w_2(x) \}$ (resp. $\mathcal{B}_\epsilon=\{w_3(x), w_4(x) \}$) does not have a limit when the parameter tends to a value for which there are logarithmic terms at the origin (resp. at $x=\epsilon$). For $\epsilon \in S^+$, there are values of $\epsilon$ for which $w_1(x)$ or $w_4(x)$ may not be defined, whereas $w_2(x)$ or $w_3(x)$ may not be defined for some values of $\epsilon$ in $S^-$.  This means that $\mathcal{B}_0$ and $\mathcal{B}_\epsilon$ are not optimal bases to describe the dynamics for all values of $\epsilon$ in the sectors $S^±$. We will rather consider the bases $\mathcal{B}^+=\{w_2(x), w_3(x) \}$ on $S^+$ and $\mathcal{B}^-=\{w_4(x), w_1(x) \}$ on $S^-$. With these bases we will explain the occurence of logarithmic terms (a phenomenon occuring for discrete values of the confluence parameter) in a continuous way. The following lemma will allow us to consider only one of the bases, namely $\mathcal{B}^+$ with $\epsilon \in S^+$. 

\begin{lemma} \label{L:e-}
The equation (\ref{E:hypergeo e}) is invariant under 
\begin{equation}
\begin{cases}
c' = 1-c+a+b \\
\epsilon'=\frac{1}{1-c'} \\
x'=\epsilon ' (1-\frac{x}{\epsilon})\\
a'=a\\
b'=b
\end{cases}
\end{equation}
which transforms $S^+$ into $S^-$ and $\mathcal{B}^+$ into $\mathcal{B}^-$.
\end{lemma}

\subsection{The confluent hypergeometric equation and its summable solutions}

Taking the limit $\epsilon \to 0$ in (\ref{E:hypergeo e}), we obtain a confluent hypergeometric equation:
\begin{equation}\label{E:conf}
x^2 \, w''(x)+\{ 1+ (1+a+b)x \} \, w'(x)+a b \, w(x)=0.
\end{equation}

A basis of solutions around the origin is 
\begin{equation}\label{E:kg}
\begin{cases}
\hat{g}(x)=\, _2F_0(a,b;-x), \\
\hat{k}(x)=e^{\frac{1}{x}}x^{1-a-b}\, _2F_0(1-a,1-b;x)=e^{\frac{1}{x}}x^{1-a-b} \hat{h}(x).
\end{cases}
\end{equation}

\begin{remark}
The confluent equation in the literature is often studied with the irregular singular point at infinity:
\begin{equation}\label{eq hyp confluente originale 2}
z u''(z)+(c'-z)u'(z)-a u(z)=0.
\end{equation}
The following transformation applied to (\ref{eq hyp confluente originale 2}) yields the confluent equation (\ref{E:conf}): 
\begin{equation}\label{E:chang w-u conf}
\begin{cases}
z=\frac{1}{x},\\
u(\frac{1}{x})=x^{a}w(x),\\
c'=a+1-b.
\end{cases}
\end{equation}
\end{remark}

\begin{figure}[h!]
\begin{center}
{\psfrag{B}{\small{$g(x)$}}
\psfrag{C}{\small{$h(x)$}}
\psfrag{A}{\footnotesize{$0$}}
\psfrag{D}{\footnotesize{$0$}}
\includegraphics[width=4cm]{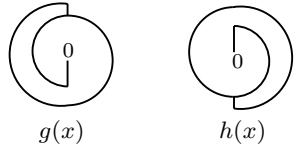}}
    \caption{Domains of the Borel sums of the confluent series $g(x)$ and $h(x)$}
    \label{fig:6}
\end{center}
\end{figure}

The following theorem is well-known, one can refer for instance to~\cite{jRjM}.
\begin{theorem} \label{T:div}
The series $\hat{g}(x)$ is divergent if and only if $a \notin -\mathbb{N}$ and $b \notin -\mathbb{N}$. It is $1$-summable in all directions except $\mathbb{R^-}$.
The series $\hat{h}(x)$ is divergent if and only if $1-a \notin -\mathbb{N}$ and $1-b \notin -\mathbb{N}$. It is $1$-summable in all directions except $\mathbb{R^+}$.
The Borel sums of these series, denoted $g(x)$ and $h(x)$, are thus defined in the sectors illustrated in Figure \ref{fig:6}.
\end{theorem}

As illustrated in Figure \ref{fig:6}, we have one Borel sum $g(x)$ in the region $\Re(x)>0$. When extending $g(x)$ to the region $\Re(x)<0$ by turning around the origin in the positive (resp. negative) direction, we get a sum $g^+(x)$ (resp. $g^-(x)$). The functions $g^+(x)$ and $g^-(x)$ are different in general and never coincide if the series is divergent. Since $g^+(x)$ and $g^-(x)$ have the same asymptotic expansion $g(x)$, their difference is a solution of (\ref{E:conf}) which is asymptotic to $0$ in the region $\Re(x)<0$, and thus 
\begin{equation}\label{E:g Stokes lambda}
g^+(xe^{2\pi i})-g^-(x)=\lambda k(x) \quad \mbox{if } \arg(x) \in (\frac{-3\pi}{2},\frac{-\pi}{2}).
\end{equation}
Similarly, we consider $h(x)$ defined in the region $\Re(x)<0$. When we extend it by turning around the origin in the positive (resp. negative) direction, we obtain the sum $h^+(x)$ (resp. $h^-(x)$). We define 
\begin{equation}
\begin{cases}
k^+(x)=e^{\frac{1}{x}}x^{1-a-b} h^+(x) \\
k^-(x)=e^{\frac{1}{x}}x^{1-a-b} h^-(x) 
\end{cases}
\end{equation}
for $\Re(x)>0$,
and
\begin{equation}
k(x)=e^{\frac{1}{x}}x^{1-a-b} h(x) 
\end{equation}
for $\Re(x)<0$.
Then we can write 
\begin{equation}\label{E:k Stokes mu}
k^+(x)-e^{2 \pi i (1-a-b)}k^-(x e^{-2\pi i})=\mu g(x) \quad \mbox{if } \arg(x) \in (\frac{-\pi}{2},\frac{\pi}{2}).
\end{equation}

\begin{remark}
For all $n \in \mathbb{Z}$, it is possible to construct a function $g_n(x)$, corresponding to the Borel sum of the divergent series $\hat{g}(x)$ in the regions $\arg(x) \in (\frac{-\pi}{2}+2\pi n,\frac{\pi}{2}+2\pi n)$. Then, $g^+_n(x)$ (resp. $g^-_n(x)$) denotes its analytic continuation in the positive (resp. negative) direction around the origin, defined in the region $\arg(x) \in (\frac{\pi}{2}+2\pi n,\frac{3\pi}{2}+2\pi n)$ (resp. $\arg(x) \in (\frac{-3\pi}{2}+2\pi n,\frac{-\pi}{2}+2\pi n)$). Since $g^+_{n+1}(xe^{2 \pi i})=g^+_n(x)$, $g^-_{n+1}(xe^{2 \pi i})=g^-_n(x)$ and $g_{n+1}(xe^{2 \pi i})=g_n(x)$, the subscript $n$ is not necessary and the functions $g(x)$, $g^+(x)$ and $g^-(x)$ are univalued. But what is important is that, when considering $g^+(x)$, the $+$ does not refer to the values of $\arg(x)$, but to the fact that $g^+(x)$ has been obtained by analytic continuation of $g(x)$ when turning in the positive direction. Similar relations for $h^+(x)$, $h^-(x)$ and $h(x)$ imply that these functions are also univalued. On the other hand, $x^{1-a-b}$ is a multivalued function, which becomes univalued as soon as $\arg(x)$ is determined.
\end{remark}

\begin{definition}
In the relations (\ref{E:g Stokes lambda}) and (\ref{E:k Stokes mu}), we call $\lambda$ and $\mu$ the Stokes multipliers associated respectively to the solutions $g(x)$ and $k(x)$. 
\end{definition}

Their values are calculated in \cite{jRjM}. Using the change of variable (\ref{E:chang w-u conf}), we have
\begin{equation}\label{E:lambda bis}
\lambda=-\frac{2 \pi i e^{i \pi (1-a-b)}}{\Gamma(a) \Gamma(b)}
\end{equation}
and
\begin{equation}\label{E:mu bis}
\mu=-\frac{2 i \pi }{\Gamma (1-a) \Gamma (1-b)}.
\end{equation}

\begin{notation}
Let us write
\begin{equation}
H^{0}(x)=\begin{cases}
\frac{k(x)}{g^-(x)} \quad  \mbox{if } \Re(x)<0\\ 
\frac{k^+(x)}{g(x)} \quad  \mbox{if } \Re(x)>0
\end{cases}
\end{equation}
and
\begin{equation}
{H^{0}}'(x)=\begin{cases}
\frac{k^-(x)}{g(x)} \quad \mbox{if }  \Re(x)>0\\
\frac{k(x)}{g^+(x)} \quad  \mbox{if }  \Re(x)<0
\end{cases}
\end{equation}
with $H^{0}(x)$ (resp. ${H^{0}}'(x)$) analytic in the complex plane minus a cut with values in $\mathbb{CP}^1$, as illustrated in Figure \ref{fig:14}. On purpose we leave the ambiguity in the argument. In this form, $H^{0}(x)$ and ${H^{0}}'(x)$ are multivalued. They will become univalued when $\arg(x)$ is specified.  
\begin{figure}[h!]
\begin{center}
{\psfrag{D}{\small{${H^{0}}'(x)$}}
\psfrag{C}{\small{$H^{0}(x)$}}
\psfrag{A}{\footnotesize{$0$}}
\psfrag{B}{\footnotesize{$0$}}
\includegraphics[width=6cm]{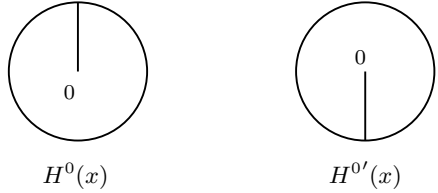}}
    \caption{Domains of $H^{0}(x)$ and ${H^{0}}'(x)$, with arbitrary radius}
    \label{fig:14}
\end{center}
\end{figure}
\end{notation}

\begin{proposition}\label{P:stokes first int}
The Stokes multiplier of $g(x)$ is
\begin{equation}\label{E:lambda}
\begin{array}{lll}
\lambda =\frac{1}{{H^{0}}'(x)}-\frac{1}{H^{0}(x)} \quad \mbox{if} \quad \arg(x) \in (\frac{-3\pi}{2},\frac{-\pi}{2}),
\end{array}
\end{equation}
while the Stokes multiplier of $k(x)$ is
\begin{equation}\label{E:mu}
\begin{array}{lll}
\mu =H^{0}(x)-e^{2 \pi i (1-a-b)}{H^{0}}'(x e^{-2 \pi i}) \quad \mbox{if} \quad \arg(x) \in (\frac{-\pi}{2},\frac{\pi}{2}).
\end{array}
\end{equation}
\end{proposition}

\begin{proof}
We have
\begin{equation}
\begin{array}{lll}
\lambda &= \frac{g^+(x e^{2 \pi i})}{k(x)}-\frac{g^-(x)}{k(x)} \\
&= \frac{g^+(x)}{k(x)}-\frac{g^-(x)}{k(x)}\\
&=\frac{1}{{H^{0}}'(x)}-\frac{1}{H^{0}(x)} \quad \mbox{if} \quad \arg(x) \in (\frac{-3\pi}{2},\frac{-\pi}{2})
\end{array}
\end{equation}
and
\begin{equation}
\begin{array}{lll}
\mu &=\frac{k^+(x)}{g(x)}-e^{2 \pi i (1-a-b)} \frac{k^-(x e^{-2 \pi i})}{g(x)}\\
&=\frac{k^+(x)}{g(x)}-e^{2 \pi i (1-a-b)}\frac{k^-(x e^{-2 \pi i})}{g(x e^{-2 \pi i})}\\
&=H^{0}(x)-e^{2 \pi i (1-a-b)}{H^{0}}'(x e^{-2 \pi i}) \quad \mbox{if} \quad \arg(x) \in (\frac{-\pi}{2},\frac{\pi}{2}).
\end{array}
\end{equation} 
\end{proof}

In view of this proposition, it will seem natural in the next section to study the monodromy of some quotient of solutions of the hypergeometric equation (\ref{E:hypergeo e}). But before, let us explore the link between divergent series in particular solutions of the confluent differential equation and analytic continuation of series appearing in solutions of the nonconfluent equation. 

\section{Divergence and Monodromy}

\subsection{Divergence and ramification: first observations}

Let us illustrate by an example the link between the divergence of a confluent series and the ramification of its unfolded series.
\begin{example}\label{exa ramif}
The series $g(x)=\, _2F_0(a,b;-x)$ is non-summable in the direction $\mathbb{R}^-$, i.e. on the left side. By continuity, when we unfold with a small $\epsilon \in \mathbb{R}$, the unfolded functions are
\begin{equation}
g^\epsilon(x)=
\begin{cases}
\, _2F_1(a,b,a+b+\frac{1}{\epsilon};1-\frac{x}{\epsilon}) &\mbox{if } \epsilon \in S^+ \\
\, _2F_1(a,b,1-\frac{1}{\epsilon};\frac{x}{\epsilon}) &\mbox{if } \epsilon \in S^-.
\end{cases}
\end{equation}
Their analytic continuations will be ramified at the left singular point and regular at the right singular point. For the special values of $\epsilon$ for which logarithmic terms may exist in the general solution at the left singular point, this will force their existence. Indeed, for these special values of $\epsilon$, the solution either has logarithmic terms or is a polynomial, in which case it cannot be ramified.
\end{example}
 
\begin{figure}[h!]
\begin{center}
{\psfrag{A}{\footnotesize{$0$}}
\psfrag{B}{\footnotesize{$0$}}
\psfrag{C}{\footnotesize{$\epsilon$}}
\psfrag{D}{\footnotesize{$0$}}
\psfrag{E}{\footnotesize{$0$}}
\psfrag{F}{\footnotesize{$\epsilon$}}
\psfrag{G}{\footnotesize{$\epsilon$}}
\psfrag{H}{\footnotesize{$0$}}
\psfrag{I}{\footnotesize{$\epsilon$}}
\psfrag{J}{\footnotesize{$0$}}
\psfrag{K}{\footnotesize{$\, _2F_0(a,b;-x)$}}
\psfrag{L}{\footnotesize{$\, _2F_1(a,b,a+b+\frac{1}{\epsilon};1-\frac{x}{\epsilon})$}}
\psfrag{M}{\footnotesize{$\, _2F_0(1-a,1-b;x)$}}
\psfrag{N}{\footnotesize{$\, _2F_1(1-a,1-b,1+\frac{1}{\epsilon};\frac{x}{\epsilon})$}}
\psfrag{O}{\footnotesize{$\, _2F_1(a,b;1-\frac{1}{\epsilon},\frac{x}{\epsilon})$}}
\psfrag{P}{\footnotesize{$\, _2F_1(1-a,1-b;2-a-b-\frac{1}{\epsilon},1-\frac{x}{\epsilon})$}}
\psfrag{Q}{$\iff$}
\psfrag{R}{$\iff$}
\psfrag{S}{or}
\psfrag{T}{or}
\psfrag{U}{$\forall \epsilon \in S^+$}
\psfrag{V}{$\forall \epsilon \in S^-$}
\psfrag{W}{$\forall \epsilon \in S^+$}
\psfrag{X}{$\forall \epsilon \in S^-$}
\includegraphics[width=14cm]{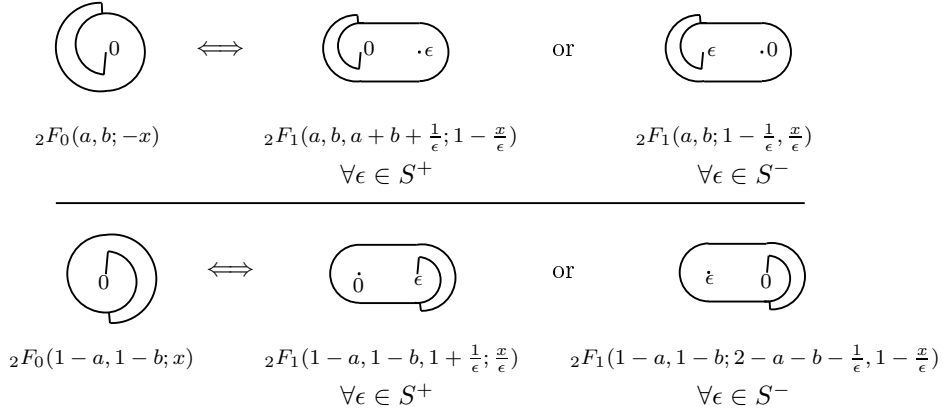}}
    \caption{Link between ramification of the analytic continuation of the hypergeometric series in the unfolded case and divergence (ramification) of the associated confluent series}
    \label{fig:1}
\end{center}
\end{figure}

This example illustrates that a direction of non-summability for a confluent series determines which merging singular point is "pathologic" (with $\epsilon$ in $S^±$) for an unfolded solution, as illustrated in Figure \ref{fig:1}.
Although subtleties are needed to adapt Example \ref{exa ramif} to the other solution $k(x)=e^{\frac{1}{x}}x^{1-a-b} h(x)$ because of the ramification of $x^{1-a-b}$, we have a similar phenomenon if we define adequately the pathology. For example, if $\epsilon \in S^+$, the singular point $x=0$ will be defined pathologic for the solution $w_3(x)$ if the analytic continuation of this solution is not an eigenvector of the monodromy operator $M_0$. This will be studied more precisely in Section \ref{S:div nondiag} using the results we will obtain in the next two sections.

\subsection{Limit of quotients of solutions on $S^±$}
\begin{figure}[h!]
\begin{center}
{\psfrag{A}{\footnotesize{$0$}}
\psfrag{B}{\footnotesize{$\epsilon$}}
\includegraphics[width=3cm]{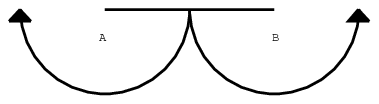}}
    \caption{Analytic continuation of $\kappa^+(\epsilon)(\frac{x}{\epsilon})^{\frac{1}{\epsilon}}(1-\frac{x}{\epsilon})^{1-\frac{1}{\epsilon}-a-b}$ for $\epsilon \in S^+$}
    \label{fig:17}
\end{center}
\end{figure}

We will later see that a divergent series in the basis of solutions at the confluence necessarily implies the presence of an obstruction that prevents an eigenvector of $M_0$ to be an eigenvector of $M_\epsilon$. As a tool for our study, we will consider the behavior of the analytic continuation of some functions of the particular solutions $w_i(x) \in \mathcal{B}^±$ when turning around singular points. A first motivation for studying these functions comes from Proposition \ref{P:stokes first int}. We will also see in Section \ref{S:ric} that these quantities have the same ramification as first integrals of a Riccati system related to the hypergeometric equation, these first integrals having a limit when $\epsilon \to 0$ on $S^±$. They are defined by
\begin{equation}
H^{\epsilon^+}(x)=\frac{\kappa^+(\epsilon) w_2(x)}{w_3(x)} \quad  \mbox{if } \epsilon \in S^+
\end{equation}
and
\begin{equation}
H^{\epsilon^-}(x)=\frac{\kappa^-(\epsilon) w_4(x)}{w_1(x)} \quad  \mbox{if } \epsilon \in S^-
\end{equation}
with  
\begin{equation}\label{E:def k}
\kappa^+(\epsilon)=\epsilon^{1-a-b} e^{\pi i(a+b-1+ \frac{1}{\epsilon})}, \quad \kappa^-(\epsilon)=\epsilon^{1-a-b} e^{-\pi i(a+b-1+ \frac{1}{\epsilon})}.
\end{equation}
$H^{\epsilon^\pm}(x)$ are first defined in $B(0,\epsilon) \cap B(\epsilon,\epsilon)$ and then analytically extended as in Figures \ref{fig:17} and \ref{fig:18}.
The coefficients $\kappa^±$ in the functions $H^{\epsilon^±}(x)$ are chosen so that $H^{\epsilon^±}(x)$ have the limit $H^{0}(x)$ when $\epsilon \to 0$ inside $S^±$. More precisely, for $\epsilon \in S^+$, we replace $f(x)=(\frac{x}{\epsilon})^{\frac{1}{\epsilon}}(1-\frac{x}{\epsilon})^{1-\frac{1}{\epsilon}-a-b}$ by  $\kappa^+(\epsilon)f(x)$, so that  the limit when $\epsilon \to 0$ and $\epsilon \in S^+$ exists and corresponds to $e^{\frac{1}{x}}x^{1-a-b}$. The limit is uniform on any simply connected compact set which does not contain $0$. The constant $\kappa^+(\epsilon)$ (resp. $\kappa^-(\epsilon)$) is the natural one to consider for $\epsilon \in S^+$ (resp. $\epsilon \in S^-$) when the analytic continuation of $\kappa^+(\epsilon)f(x)$ (resp. $\kappa^-(\epsilon)f(x)$) is done like in Figure \ref{fig:17} (resp. Figure \ref{fig:18}).

\begin{figure}[h!]
\begin{center}
{\psfrag{A}{\footnotesize{$\epsilon$}}
\psfrag{B}{\footnotesize{$0$}}
\includegraphics[width=3cm]{FigD}}
    \caption{Analytic continuation of $\kappa^-(\epsilon)(\frac{x}{\epsilon})^{\frac{1}{\epsilon}}(1-\frac{x}{\epsilon})^{1-\frac{1}{\epsilon}-a-b}$ for $\epsilon \in S^-$}
    \label{fig:18}
\end{center}
\end{figure}

\begin{proposition}\label{H limit prop}
When $\epsilon \to 0$ and $\epsilon \in S^+$ (resp. $\epsilon \in S^-$), $H^{\epsilon^+}(x)$ (resp. $H^{\epsilon^-}(x)$) converges uniformly to $H^{0}(x)$ on any simply connected compact subset of the domain of $H^{0}(x)$ illustrated in Figure \ref{fig:14}. More precisely, we have the uniform limits on compact subsets:
\begin{equation}\label{E:conv unif results}
\begin{array}{lll}
\begin{cases}
\lim_{\begin{subarray}{l}\epsilon \to 0\\ \epsilon \in S^+\end{subarray}}\kappa^+(\epsilon) w_2(x)=k^+(x)\\
\lim_{\begin{subarray}{l}\epsilon \to 0\\ \epsilon \in S^+\end{subarray}}w_3(x)=g(x)
\end{cases}
\begin{cases}
\lim_{\begin{subarray}{l}\epsilon \to 0\\ \epsilon \in S^-\end{subarray}}\kappa^-(\epsilon) w_4(x)=k^+(x)\\
\lim_{\begin{subarray}{l}\epsilon \to 0\\ \epsilon \in S^-\end{subarray}}w_1(x)=g(x)
\end{cases}
\end{array}
\end{equation}
\end{proposition}

\begin{proof}
The hypergeometric functions appearing in $w_k(x)$ ($k=1,2,3,4$) and having the limit $h(x)$ or $g(x)$ are ramified as illustrated in Figure \ref{fig:1}, which suggests to take sectors like in Figure \ref{fig:14} when considering the quotient of these functions. 

We first prove the uniform convergence $w_3(x)$ to $g(x)$ on simply connected compact subsets of the domain $\{x, |\arg(x)|<\frac{3\pi}{2}\}$ for $\epsilon \in S^+$. This proof has been inspired by~\cite{cZ94}.
Let us suppose that $a-b \notin \mathbb{Z}$. The Borel sum of $g(x)$ is the same as the analytic continuation of this solution, which is (see~\cite{yL69})

\begin{equation}
w_3(x)=\frac{\Gamma(a+b+\frac{1}{\epsilon})\Gamma(b-a)}{\Gamma(b) \Gamma(b+\frac{1}{\epsilon})}w_5(x) +\frac{\Gamma(a+b+\frac{1}{\epsilon})\Gamma(a-b)}{\Gamma(a) \Gamma(a+\frac{1}{\epsilon})}w_6(x)
\end{equation}
with 
\begin{equation}
\begin{cases}
w_5(x)=(\frac{\epsilon}{x})^a \,_2F_1(a,a+\frac{1}{\epsilon}, a+1-b; \frac{\epsilon}{x}) \\
w_6(x)=(\frac{\epsilon}{x})^b \,_2F_1(b,b+\frac{1}{\epsilon}, b+1-a; \frac{\epsilon}{x}).
\end{cases}
\end{equation}

The function $_2F_1(a,a+\frac{1}{\epsilon}, a+1-b; \frac{\epsilon}{x})$ converges uniformly on simply connected compact subsets to $\,_1F_1(a,a+1-b; \frac{1}{x})$ and we have
\begin{equation}
\begin{array}{lll}
\lim_{\begin{subarray}{l}\epsilon \to 0 \\ \epsilon \in S^+ \end{subarray}}\frac{\epsilon^a \Gamma(a+b+\frac{1}{\epsilon})}{\Gamma(b+\frac{1}{\epsilon})}=1.
\end{array}
\end{equation}
The same relations apply with $a$ and $b$ interchanged so $w_3(x)$ converges uniformly on simply connected compact subsets to
\begin{equation}
g(x)=\frac{\Gamma(b-a)}{\Gamma(b) }x^{-a} \,_1F_1(a, a+1-b; \frac{1}{x}) +\frac{\Gamma(a-b)}{\Gamma(a)}x^{-b} \,_1F_1(b, b+1-a; \frac{1}{x})
\end{equation}

Let us suppose now that $a-b=-m$ with $m \in \mathbb{N}$. We take $h$ small, we let $a=b-m+h$. We first show that $\lim_{h \to 0}w_3(x)$ exists with $x$ on a simply connected compact subset of the domain $\{x, |\arg(x)|<\frac{3\pi}{2}\}$. We write $w_3(x)$ as
\begin{equation}
\begin{array}{llll}\label{E:w3 de h}
w_3(x)=(a-b)\Gamma(b-a)\Gamma(a-b)&\Gamma(a+b+\frac{1}{\epsilon}) \\&\left[ \frac{w_5(x)}{\Gamma(b)\Gamma(b+\frac{1}{\epsilon}) \Gamma(a-b+1)} - \frac{w_6(x)}{\Gamma(a)\Gamma(a+\frac{1}{\epsilon}) \Gamma(b-a+1)}\right]
\end{array}
\end{equation}
and take the limit $h \to 0$ with $a=b-m+h$. The part inside brackets has a zero at $h=0$ since
\begin{equation}
\begin{array}{lll}
\lim_{h \to 0}\frac{w_5(x)}{\Gamma(a-b+1)}&=(\frac{\epsilon}{x})^{b}\frac{(b-m)_m (b-m+\frac{1}{\epsilon})_m}{m!}\, _2F_1(b, b+\frac{1}{\epsilon},m+1;\frac{\epsilon}{x})\\
&=\frac{\Gamma(b) \Gamma(b+\frac{1}{\epsilon})w_6(x)}{\Gamma(a)\Gamma(a+\frac{1}{\epsilon})\Gamma(b-a+1)}
\end{array}
\end{equation}
The left part of (\ref{E:w3 de h}) has a simple pole at $h=0$ so $\lim_{h \to 0}w_3(x)$ exists. Since $w_3(x)$ is an analytic function of $h$ on a punctured neighborhood of $h=0$, we have that $w_3(x)$ converges uniformly on simply connected compact subsets to $\lim_{h \to 0}w_3(x)$ when $h \to 0$. Similarly, $g(x)$ converges uniformly on simply connected compact subsets to $\lim_{h \to 0}g(x)$ since 
\begin{equation}
\begin{array}{lll}
\lim_{h \to 0}\frac{_1F_1(a,a+1-b;\frac{1}{x})}{x^a \Gamma(a-b+1)}=\frac{\Gamma(b)_1F_1(b,b+1-a;\frac{1}{x})}{x^b\Gamma(a) \Gamma(b-a+1)}
.\end{array}
\end{equation}
Hence,  $\lim_{h \to 0}w_3(x)$ converges uniformly on simply connected compact subsets to $\lim_{h \to 0}g(x)$ when $\epsilon \to 0$ with $\epsilon \in S^+$.  Interchanging $a$ and $b$ leads to the case $b-a \in -\mathbb{N}$. 

Now, $w_2(x)$ (as in (\ref{E:w1(x) w2(x)})) converges uniformly to $k(x)$ on simply connected compact subsets of the domain $\{x, |\arg(-x)|<\frac{3\pi}{2}\}$ to $k(x)$. Indeed, we can decompose $\kappa^+(\epsilon) w_2(x)$ as
\begin{equation}
\left( e^{\frac{\pi i}{\epsilon}}(\frac{x}{\epsilon})^{\frac{1}{\epsilon}} (1-\frac{x}{\epsilon})^{-\frac{1}{\epsilon}}\right) \left( (x-\epsilon)^{1-a-b}\, _2F_1(1-a,1-b,1+\frac{1}{\epsilon};\frac{x}{\epsilon}) \right).
\end{equation}

The first part converges to $e^{\frac{1}{x}}$. The second part converges to $x^{1-a-b} \,_2F_0(1-a,1-b;x)$. The fact that $_2F_1(1-a,1-b,1+\frac{1}{\epsilon};\frac{x}{\epsilon})$ converges uniformly on simply connected compact subsets to $_2F_0(1-a,1-b;x)$ can be obtained from the convergence of $w_3(x)$ to $g(x)$ by a change of coordinates. The case $\epsilon \in S^-$ is similar. 
\end{proof}

\subsection{Divergence and nondiagonal form of the monodromy operator in the basis $\mathcal{B}^+$}\label{S:div nondiag}

It is clear that $w_2(x)$ is an eigenvector of the monodromy operator $M_0$ with eigenvalue $e^{\frac{i \pi}{\epsilon}}$, and that $w_3(x)$ is an eigenvector of $M_\epsilon$ with eigenvalue $1$. In general, eigenvectors of the monodromy operators $M_0$ and $M_\epsilon$ should not coincide. In the generic case, the analytic continuation of an eigenvector of the monodromy operator $M_0$ is not an eigenvector of $M_\epsilon$. If we are in the generic case and this persists to the limit $\epsilon=0$, then at the limit we have a nonzero Stokes multiplier. The results stated in the next theorem tell us whether or not the analytic continuation of $w_3(x)$ (resp. $w_2(x)$) is an eigenvector of $M_0$ (resp. $M_\epsilon$). This is done in the two covering sectors $S^±$ of a small neighborhood of $\epsilon$, and it includes the presence of logarithmic terms: we will detail this last part in Theorem \ref{T:logarithmic terms} below.

\begin{figure}[h!]
\begin{center}
{\psfrag{A}{\footnotesize{$0$}}
\psfrag{B}{\footnotesize{$\epsilon$}}
\psfrag{C}{\small{$w_{(0, \pi)}(x)$}}
\psfrag{D}{\small{$w_{(0, -\pi)}(x)$}}
\psfrag{E}{\small{$w_{(\epsilon, -\pi)}(x)$}}
\psfrag{F}{\small{$w_{(\epsilon, \pi)}(x)$}}
\includegraphics[width=7cm]{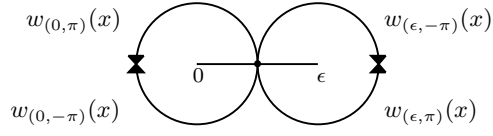}}
    \caption{Analytic continuation of $w(x)$}
    \label{fig:3}
\end{center}
\end{figure}

\begin{notation}\label{N:analytic ext}
Let $w_{(\delta, \theta)} (x)$ be the analytic continuation of $w(x)$ when starting on $(0,\epsilon)$ and turning of an angle $\theta$ around $x=\delta$, with $\delta \in \{0, \epsilon \}$ (see Figure \ref{fig:3}). In short, $w_{(\delta, \pi)} (x)$ can be obtained from the action of the monodromy operator around $x=\delta$ applied on $w_{(\delta, -\pi)} (x)$.
\end{notation}

\begin{theorem} \label{Stokes coef limit w2 w3}
\begin{itemize}
 \item If $\epsilon \in S^+$, then
\begin{equation}\label{E:monodromymatrix 1}
\begin{pmatrix}
\kappa^+(\epsilon)w_{2,(0,\pi)}\\
w_{3,(0,\pi)}
\end{pmatrix}=
\begin{pmatrix}
e^{\frac{2\pi i}{\epsilon}} &0 \\
\lambda^+(\epsilon) &1
\end{pmatrix}
\begin{pmatrix}
\kappa^+(\epsilon)w_{2,(0,-\pi)}\\
w_{3,(0,-\pi)}
\end{pmatrix}
\end{equation}
and
\begin{equation}\label{E:monodromymatrix 2}
\begin{pmatrix}
\kappa^+(\epsilon)w_{2,(\epsilon ,\pi)}\\
w_{3,(\epsilon ,\pi)}
\end{pmatrix}=
\begin{pmatrix}
e^{2\pi i(1-a-b-\frac{1}{\epsilon})} &\mu^+(\epsilon) \\
0 &1
\end{pmatrix}
\begin{pmatrix}
\kappa^+(\epsilon) w_{2,(\epsilon ,-\pi)}\\
w_{3,(\epsilon ,-\pi)} 
\end{pmatrix},
\end{equation}
with
\begin{equation}\label{E: mu plus}
\mu^+(\epsilon )=\frac{-2 \pi i}{\Gamma (1-a) \Gamma (1-b)}  \frac{\epsilon ^{1-a-b} \Gamma (1+ \frac{1}{\epsilon})}{ \Gamma (a+b+\frac{1}{\epsilon})} 
\end{equation}
and
\begin{equation}\label{E: lambda plus}
\lambda^+(\epsilon )=\frac{-2 \pi i e^{\pi i (1-a-b)}}{\Gamma (a) \Gamma (b) }  \frac{\epsilon ^{a+b-1} \Gamma (a+b+\frac{1}{\epsilon})}{\Gamma (1+ \frac{1}{\epsilon})}.
\end{equation}
Hence, when it is nonzero, the coefficient $\lambda^+(\epsilon)$ (resp. $\mu^+(\epsilon)$) represents the obstruction that prevents $w_3(x)$ (resp. $w_2(x)$) of being an eigenvector of the monodromy operator around $x=0$ (resp. $x=\epsilon$).

\item If $\epsilon \in S^-$, then
\begin{equation}\label{E:monodromymatrix 1 -}
\begin{pmatrix}
\kappa^-(\epsilon) w_{4,(\epsilon,\pi)}\\
w_{1,(\epsilon,\pi)}
\end{pmatrix}=
\begin{pmatrix}
e^{2\pi i (1-\frac{1}{\epsilon}-a-b)} &0 \\
\lambda^-(\epsilon) &1
\end{pmatrix}
\begin{pmatrix}
\kappa^-(\epsilon) w_{4,(\epsilon,-\pi)}\\
w_{1,(\epsilon,-\pi)}
\end{pmatrix}
\end{equation}
and
\begin{equation}\label{E:monodromymatrix 2 -}
\begin{pmatrix}
\kappa^-(\epsilon)w_{4,(0,\pi)}\\
w_{1,(0,\pi)}
\end{pmatrix}=
\begin{pmatrix}
e^{\frac{2\pi i}{\epsilon}} &\mu^-(\epsilon) \\
0 &1
\end{pmatrix}
\begin{pmatrix}
\kappa^-(\epsilon) w_{4,(0,-\pi)}\\
w_{1,(0,-\pi)} 
\end{pmatrix},
\end{equation}
with 
\begin{equation}\label{E: mu moins}
\mu^-(\epsilon )=\frac{-2 \pi i}{\Gamma (1-a) \Gamma (1-b)}  \frac{(\epsilon e^{\pi i}) ^{1-a-b} \Gamma (2- \frac{1}{\epsilon}-a-b)}{ \Gamma (1-\frac{1}{\epsilon})} 
\end{equation}
and
\begin{equation}\label{E: lambda moins}
\lambda^-(\epsilon )=\frac{-2 \pi i}{\Gamma (a) \Gamma (b) }  \frac{(\epsilon ) ^{a+b-1} \Gamma (1-\frac{1}{\epsilon})}{\Gamma (2-\frac{1}{\epsilon}-a-b)}  .
\end{equation}

Hence, when it is nonzero, the coefficient $\lambda^-(\epsilon)$ (resp. $\mu^-(\epsilon)$) represents the obstruction that prevents $w_1(x)$ (resp. $w_4(x)$) of being an eigenvector of the monodromy operator around $x=\epsilon$ (resp. $x=0$).
\end{itemize}

Then, with the limit taken for any path in $S^+$ or in $S^-$, we have
\begin{equation}
\lim_{\epsilon \to 0}\mu^±(\epsilon )=\mu
\end{equation}
and 
\begin{equation}
\lim_{\epsilon \to 0}\lambda^±(\epsilon )=\lambda,
\end{equation}
which are precisely the Stokes multipliers associated to the solutions $k(x)$ and $g(x)$ and given by (\ref{E:lambda bis}) and (\ref{E:mu bis}).
\end{theorem}

\begin{proof}
Let $\epsilon \in S^+$. To make analytic continuation of the solutions $w_2(x)$ and $w_3(x)$, we need to make further restrictions on the values of $\epsilon$, but we will shortly show the validity of the result without these hypotheses. We have (see for example \cite{yL69})
\begin{itemize}
\item
if $2-\frac{1}{\epsilon}-a-b \notin -\mathbb{N}$,
\begin{equation}\label{prolong w2}
\begin{array}{lllll}
w_2(x)&=\frac{\Gamma (1-\frac{1}{\epsilon}-a-b) \Gamma (1+\frac{1}{\epsilon})}{\Gamma (1-a) \Gamma (1-b)} &w_3(x) + \frac{\Gamma (a+b-1+\frac{1}{\epsilon}) \Gamma (1+\frac{1}{\epsilon})}{\Gamma (a+\frac{1}{\epsilon}) \Gamma (b+\frac{1}{\epsilon})} &w_4(x)\\
&=D(\epsilon) &w_3(x) + E(\epsilon) &w_4(x) ;
\end{array}
\end{equation}
\item if $1-\frac{1}{\epsilon} \notin -\mathbb{N}$,
\begin{equation}
\begin{array}{lllll}
w_3(x)&=\frac{\Gamma (\frac{1}{\epsilon}) \Gamma (a+b+\frac{1}{\epsilon})}{\Gamma (b+\frac{1}{\epsilon}) \Gamma (a+\frac{1}{\epsilon})} &w_1(x) + \frac{\Gamma (a+b+\frac{1}{\epsilon}) \Gamma (-\frac{1}{\epsilon})}{\Gamma (a) \Gamma (b)} &w_2(x)\\
&=A(\epsilon) &w_1(x) + B(\epsilon) &w_2(x) .
\end{array}
\end{equation}
\end{itemize}

These relations allow the calculation of the monodromy of $w_2(x)$ (resp. $w_3(x)$) around $x=\epsilon$ (resp. $x=0$). The explosion of the coefficients (coefficients becoming infinite) for specific values of $\epsilon$ corresponds to the presence of logarithmic terms in the general solution around the singular point $x=\epsilon$ (resp. $x=0$). 
We have, in the region $B(0,\epsilon) \cap B(\epsilon,\epsilon)$ (with the hypothesis that $2-\frac{1}{\epsilon}-a-b \notin -\mathbb{N}$),
\begin{equation}
\begin{array}{lll}
\kappa^+(\epsilon) w_2(x)&=&\kappa^+(\epsilon) ( D(\epsilon) w_3(x) + E(\epsilon) w_4(x))\\
  &=&\kappa^+(\epsilon) ( D(\epsilon ) \, _2F_1(a,b,a+b+\frac{1}{\epsilon};1-\frac{x}{\epsilon}) \\
&&+ E(\epsilon )  (\frac{x}{\epsilon})^{\frac{1}{\epsilon}}(1-\frac{x}{\epsilon})^{1-\frac{1}{\epsilon}-a-b} \, _2F_1(1-a,1-b,-\frac{1}{\epsilon}+2-a-b;1-\frac{x}{\epsilon})).
\end{array}
\end{equation}
Since $w_{3,(\epsilon,-\pi)}=w_{3,(\epsilon,\pi)}$, we obtain
\begin{equation}
\kappa^+(\epsilon) w_{2,(\epsilon,\pi)}=e^{2 \pi i (1-a-b-\frac{1}{\epsilon})}\kappa^+(\epsilon) w_{2,(\epsilon,-\pi)}+\mu^+(\epsilon )w_{3,(\epsilon,-\pi)}
\end{equation}
with
\begin{equation}
\begin{array}{lll}
\mu^+(\epsilon)&=D(\epsilon )\epsilon ^{1-a-b} e^{ \pi i (a+b-1+\frac{1}{\epsilon})} \left( 1-e^{ 2 \pi i (1-a-b-\frac{1}{\epsilon})} \right)\\
&=-D(\epsilon )\epsilon ^{1-a-b} \left( e^{ \pi i (1-a-b-\frac{1}{\epsilon})}-e^{ - \pi i (1-a-b-\frac{1}{\epsilon})} \right).
\end{array}
\end{equation}

Since $\sin(z)=\frac{e^{i z}-e^{-iz}}{2i}$ and $\Gamma(z)\sin(\pi z) =\frac{\pi}{\Gamma (1-z)}$, we can simplify the latter expression:
\begin{equation}
\begin{array}{lll}
\mu^+(\epsilon )&=-2 i D(\epsilon ) \epsilon ^{1-a-b} \sin(\pi (1-a-b-\frac{1}{\epsilon})\\
&=-2 i \frac{\Gamma (1- \frac{1}{\epsilon}-a-b) \Gamma (1+ \frac{1}{\epsilon})}{\Gamma (1-a) \Gamma (1-b)}\epsilon ^{1-a-b}  \sin (\pi (1-a-b-\frac{1}{\epsilon})) \\
&=-2 \pi i \frac{\Gamma (1+ \frac{1}{\epsilon})}{\Gamma (1-a) \Gamma (1-b)} \epsilon ^{1-a-b}\frac{1}{\Gamma (a+b+\frac{1}{\epsilon})} .
\end{array}
\end{equation}

Remark that this expression is defined even if $2-\frac{1}{\epsilon}-a-b \in -\mathbb{N}$, so we have removed the indeterminacy!

In the particular case $a+b \in \mathbb{Z}$,
\begin{equation}
\mu^+(\epsilon )=-\frac{-2 i \pi }{\Gamma (1-a) \Gamma (1-b)} \epsilon ^{1-a-b} r(a+b)
\end{equation}
with
\begin{equation}
\begin{array}{lll}
r(\gamma)&=\frac{\Gamma (1+ \frac{1}{\epsilon})}{\Gamma (\gamma +\frac{1}{\epsilon})}
&=
\begin{cases}
\prod_{j=1} ^{\gamma-1} \frac{1}{\frac{1}{\epsilon}+j}  & \gamma>1 \quad , \\
\prod_{j=\gamma} ^{0} (\frac{1}{\epsilon}+j) & \gamma<1 \quad , \\
1 & \gamma=1 \quad .
\end{cases}
\end{array}
\end{equation}

Finally,
\begin{equation}
\lim_{\begin{subarray}{l}\epsilon \to 0\\ \epsilon \in S^+\end{subarray}} \epsilon^{1-a-b}\frac{\Gamma(\frac{1}{\epsilon}+1)}{\Gamma(\frac{1}{\epsilon}+a+b)}=1 .
\end{equation}
Hence
\begin{equation}
\lim_{\begin{subarray}{l}\epsilon \to 0\\ \epsilon \in S^+\end{subarray}}\mu^+(\epsilon )=-\frac{2 i \pi }{\Gamma (1-a) \Gamma (1-b)}=\mu.
\end{equation}

Let $\epsilon_n$ such that $2-\frac{1}{\epsilon_n}-a-b =-n$, $n \in \mathbb{N}$. Recall that we have supposed $\epsilon \ne \epsilon_n$ to obtain $\mu^+(\epsilon)$. Since $\mu^+(\epsilon)$ is analytic in a punctured disk $B(\epsilon_n,\rho) \backslash \{ \epsilon_n \}$ (for some well chosen $\rho \in \mathbb{R}_+$), and $\lim_{\epsilon \to \epsilon_n} \mu^+(\epsilon)$ exists, then $\mu^+(\epsilon)$ is analytic in $B(\epsilon_n,\rho)$. Hence, the result obtained is valid without the restriction $2-\frac{1}{\epsilon}-a-b \notin -\mathbb{N}$.

A similar calculation gives, with $w_{2,(0,\pi)}=e^{\frac{2 \pi i}{\epsilon}}w_{2,(0,-\pi)}$,

\begin{equation}
w_{3,(0,\pi)}=w_{3,(0,-\pi)}+\lambda^+(\epsilon )\kappa^+(\epsilon) w_{2,(0,-\pi)}
\end{equation}
with $\lambda^+(\epsilon)=B(\epsilon)e^{-\pi i (a+b-1+\frac{1}{\epsilon})}\epsilon ^{a+b-1} \left( e^{\frac{2 \pi i }{\epsilon}}-1  \right)$.

And then
\begin{equation}
\lambda^+(\epsilon )=-2 \pi i e^{\pi i (1-a-b)} \frac{1}{\Gamma (a) \Gamma (b)} \epsilon ^{a+b-1}\frac{\Gamma (a+b+\frac{1}{\epsilon})}{\Gamma (1+ \frac{1}{\epsilon})},
\end{equation}
which, for $a+b \in \mathbb{Z}$, yields
\begin{equation}
\lambda^+(\epsilon )= \frac{-2 \pi i e^{\pi i (1-a-b)} \epsilon ^{a+b-1} }{\Gamma (a) \Gamma (b)} \frac{1}{r(a+b)}.
\end{equation}

Hence,
\begin{equation}
\lim_{\begin{subarray}{l}\epsilon \to 0\\ \epsilon \in S^+\end{subarray}}\lambda^+(\epsilon )= \frac{-2 \pi i e^{i \pi (1-a-b)}}{\Gamma(a) \Gamma(b)}=\lambda.
\end{equation}

Finally, Lemma \ref{L:e-} and equation (\ref{E:c}) relates the case $\epsilon' \in S^+$ to the case $\epsilon \in S^-$, and we have, denoting $w_i(x)$ by $w_i(x,\epsilon)$, 
\begin{equation}
\begin{array}{lll}
\kappa^+(\epsilon)= ( e^{\pi i} \frac{\epsilon'}{\epsilon} )^{a+b-1} \kappa^-(\epsilon')\\
w_2(x,\epsilon)=w_4(x',\epsilon')\\
w_3(x,\epsilon)=w_1(x',\epsilon')
\end{array}
\end{equation}
\end{proof}

\begin{theorem}\label{T:logarithmic terms}
\begin{enumerate}
\item \label{premier} If the series $g(x)$ is divergent, then, for all $\epsilon \in S^+$ (resp. for all $\epsilon \in S^-$), $w_3(x)$ (resp. $w_1(x)$) is not an eigenvector of the monodromy operator $M_0$ (resp. $M_\epsilon$). In particular, this forces the existence of logarithmic terms at $x=0$ (resp. $x=\epsilon$) for all special values of $\epsilon$ for which they may exist.
\item \label{deuxieme} Conversely, for fixed $a$ and $b$, if $w_3(x)$ (resp. $w_1(x)$) is not an eigenvector of the monodromy operator $M_0$ (resp. $M_\epsilon$) for some $\epsilon \in S^+$ (resp. for some $\epsilon \in S^-$), then the series $g(x)$ is divergent. 
\item If the series $h(x)$ is divergent, then, for all $\epsilon \in S^+$ (resp. for all $\epsilon \in S^-$), $w_2(x)$ (resp. $w_4(x)$) is not an eigenvector of the monodromy operator $M_\epsilon$ (resp. $M_0$). In particular, this forces the existence of logarithmic terms at $x=\epsilon$ (resp. $x=0$) for all special values of $\epsilon$ for which they may exist.
\item Conversely, for fixed $a$ and $b$, if $w_2(x)$ (resp. $w_4(x)$) is not an eigenvector of the monodromy operator $M_\epsilon$ (resp. $M_0$) for some $\epsilon \in S^+$ (resp. for some $\epsilon \in S^-$), then the series $h(x)$ is divergent. 
\end{enumerate}
\end{theorem}

\begin{proof}
Let $\epsilon \in S^+$ (the proof for $\epsilon \in S^-$ is similar). With Theorem \ref{T:div}, we have that $g(x)$ is divergent if and only if $\lambda \ne 0$. Since $\lim_{\epsilon \to 0}\lambda^+(\epsilon)=\lambda$, we have $\lambda^+(\epsilon) \ne 0$ for $\epsilon \in S^+$ provided the radius of $S^+$ is sufficiently small. If $w_3(x)$ were an eigenvector of the monodromy operator $M_0$, then we would have $\lambda^+(\epsilon)=0$ which is a contradiction. If $\lambda^+(\epsilon) \ne 0$, then the analytic continuation of $w_3(x)$ is ramified around $x=0$. When $1-\frac{1}{\epsilon} \in -\mathbb{N}$, $w_2(x)$ is not ramified around $x=0$ and either $w_1(x)$ is a polynomial or it has logarithmic terms. Since the analytic continuation of $w_3(x)$ is ramified at $x=0$ and since it is a linear combination of $w_1(x)$ and $w_2(x)$, we are forced to have $w_1(x)$ with logarithmic terms. The argument is similar for $w_2(x)$.

To prove the converse, we use the expressions (\ref{E: mu plus}) and (\ref{E: lambda plus}): for $\epsilon \in S^+$ and $a$ and $b$ fixed, we have $\lambda^+(\epsilon) \ne 0$ if and only if $\lambda \ne 0$ as well as $\mu^+(\epsilon) \ne 0$ if and only if $\mu \ne 0$.
\end{proof}

Hence, the singular direction $\mathbb{R}^-$ (resp. $\mathbb{R}^+$) of the 1-summable series $g(x)$ (resp. $h(x)$) is directly related to the presence of logarithmic terms at the \emph{left} (resp. \emph{right}) singular point for specific values of the confluence parameter.

\begin{remark} The necessary condition (\ref{premier}) in Theorem \ref{T:logarithmic terms} is still valid when $a$ and $b$ are analytic functions $a(\epsilon)$ and $b(\epsilon)$. A counter example to the converse (\ref{deuxieme}), for instance with $a(\epsilon)$ and $b(\epsilon)$ non constant, is given by
\begin{equation}
\begin{cases}
a(\epsilon)=n+\epsilon, \quad n \in -\mathbb{N} \\
b(\epsilon)=m+\epsilon, \quad m \in \mathbb{N}^*.
\end{cases}
\end{equation}
\end{remark}

Looking at Theorem \ref{Stokes coef limit w2 w3}, it is clear that, even in the convergent case, there is some wild behavior ($e^{\frac{2 \pi i}{\epsilon}}$) in the monodromy of the solutions which does not go to the limit.  Fortunately, this wild behavior is linear. In the next section, we will separate it from the non linear part in order to get a limit for the latter.

\subsection{The wild and continous part of the monodromy operator}\label{S:separ}

In this section, we see that the monodromy of $H^{\epsilon^±}(x)$ can be separated in a wild part and continuous part. This is the advantage of studying the monodromy of $H^{\epsilon^±}(x)$ instead of the monodromy of each solution. The wild part is present even in the case of convergence of the confluent series in $g(x)$ and in $k(x)$ and is purely linear. The continous part leads us to the Stokes coefficients. This is still done in the two covering sectors $S^±$ of a small neighborhood of $\epsilon$. 

\begin{theorem} \label{Stokes coef limit} 
Let $H^{\epsilon^±}_{i,(\delta, \theta)}(x)$ be obtained from analytic continuation of $H^{\epsilon^±}(x)$ as in notation \ref{N:analytic ext}. The relation between $H^{\epsilon^±}_{(\epsilon, \mp \pi)}$ and $H^{\epsilon^±}_{(\epsilon,\pm \pi)}$, as well as the relation between $H^{\epsilon^±}_{(0,\mp \pi)}$ and $H^{\epsilon^±}_{(0,\pm \pi)}$ may be separated into
\begin{itemize}
\item a wild linear part with no limit at $\epsilon=0$
\item a continuous non linear part
\end{itemize}
on each of the sectors $S^±$.
More precisely,
\begin{itemize}
\item if $\epsilon \in S^+$, 
\begin{equation}\label{E:H tour e}
H^{\epsilon^+}_{(\epsilon,-\pi)}=e^{2 \pi i (a+b-1+\frac{1}{\epsilon})}(H^{\epsilon^+}_{(\epsilon,\pi)}-\mu^+(\epsilon ))
\end{equation}
and
\begin{equation}\label{E:H tour 0}
\frac{1}{H^{\epsilon^+}_{(0,\pi)}}=e^{\frac{-2 \pi i }{\epsilon}} \left( \frac{1}{H^{\epsilon^+}_{(0,-\pi)}}+\lambda^+(\epsilon ) \right)
\end{equation}
with $\mu^+(\epsilon )$ and $\lambda^+(\epsilon )$ as in (\ref{E: mu plus}) and (\ref{E: lambda plus}).

\item if $\epsilon \in S^-$,
\begin{equation}\label{E:K tour e}
H^{\epsilon^-}_{(0,-\pi)}=e^{\frac{-2 \pi i }{\epsilon}}(H^{\epsilon^-}_{(0,\pi)}-\mu^-(\epsilon ))
\end{equation}
and
\begin{equation}\label{E:K tour 0}
\frac{1}{H^{\epsilon^-}_{(\epsilon,\pi)}}=e^{2 \pi i (a+b-1+\frac{1}{\epsilon})} \left( \frac{1}{H^{\epsilon^-}_{(\epsilon,-\pi)}}+\lambda^-(\epsilon ) \right)
\end{equation}
with $\mu^-(\epsilon )$ and $\lambda^-(\epsilon )$ as in (\ref{E: mu moins}) and (\ref{E: lambda moins}).
\end{itemize}

\end{theorem}

\begin{proof}
The proof is a mere calculation using (\ref{E:monodromymatrix 1}), (\ref{E:monodromymatrix 2}), (\ref{E:monodromymatrix 1 -}) and (\ref{E:monodromymatrix 2 -}).
\end{proof}

\begin{proposition}
To know which invariants are realisable, it is sufficient to look at the product $\lambda^+(\epsilon)\mu^+(\epsilon)$. If $a$ and $b$ are analytic functions of $\epsilon$, this last product is analytic in a neighborhood of $\epsilon=0$.
\end{proposition}
\begin{proof}
If $\mu^+(\epsilon) \ne 0$, we can take $\mu^+(\epsilon) w_3(x)$ instead of $w_3(x)$ in the expression for $H^{\epsilon^+}(x)$. Then, $\mu^+(\epsilon)$ is replaced by $1$ in equation (\ref{E:H tour e}) and $\lambda^+(\epsilon)$ is replaced by $\lambda^+(\epsilon)\mu^+(\epsilon)$ in equation (\ref{E:H tour 0}). Similarly if $\lambda^+(\epsilon)\ne 0$. So we can regard our invariants as $1$ and $\lambda^+(\epsilon)\mu^+(\epsilon)$, instead of $\lambda^+(\epsilon)$ and $\mu^+(\epsilon)$ in the case where one of them is different from $0$.  We have
\begin{equation}
\begin{array}{lll}
\lambda^+(\epsilon)\mu^+(\epsilon)&=-\frac{4 \pi^2 e^{\pi i (1-a-b)}}{\Gamma (1-a) \Gamma (1-b) \Gamma (a) \Gamma (b)}\\
&= - 4 e^{\pi i (1-a-b)} \sin(\pi a) \sin(\pi b)\\
&= -(1-e^{-2\pi i a})(1-e^{-2\pi i b})\\
&=\lambda^-(\epsilon)\mu^-(\epsilon).
\end{array}
\end{equation}
\end{proof}

\begin{remark}
If $\mu^+(\epsilon) \ne 0$ (resp. $\lambda^+(\epsilon)\ne 0$), the product $\lambda^+(\epsilon)\mu^+(\epsilon)=\lambda^-(\epsilon)\mu^-(\epsilon)$ is zero precisely when $a \in -\mathbb{N}$  or $b \in -\mathbb{N}$ (resp. $1-a \in -\mathbb{N}$  or $1-b \in -\mathbb{N}$), i.e. when $g(x)$ (resp. $k(x)$) is a convergent solution.
\end{remark}

\begin{remark}
When $a+b=1$, we have $\mu^+(\epsilon)=\lambda^+(\epsilon)$ and $\mu^-(\epsilon)=\lambda^-(\epsilon)$ (and $\mu=\lambda$). We will see in Remark \ref{R:formal} of Section \ref{S:ric} that this is the particular case when the formal invariants of the two saddle-nodes of the Riccati equation (\ref{E:ric 1}) vanish.
\end{remark}

\section{A related Riccati system}\label{S:ric}

\subsection{First integrals of a Riccati system related to the hypergeometric equation (\ref{E:hypergeo e})}
We studied the monodromy of $H^{\epsilon^±}(x)=\frac{\kappa^±(\epsilon) w_i(x)}{w_j(x)}$ (with $(i,j)=\begin{cases}(2,3), \, \epsilon \in S^+\\
(4,1), \, \epsilon \in S^-
\end{cases}$) instead of the monodromy of each solution $w_k(x)$, for $k=i,j$. To justify this choice, we transform the hypergeometric equation into a Riccati equation (see for instance~\cite{Hille}) and find a first integral of the Riccati system.

\begin{proposition}
The Riccati system
\begin{equation}\label{E:ric 1}
\begin{cases}
\dot x=x(x-\epsilon)\\
\dot y =a b x (x-\epsilon)+(-1+(1-a-b)x)y+y^2
\end{cases}
\end{equation}
is related to the hypergeometric equation (\ref{E:hypergeo e}) with singular points at $\{0,\epsilon,\infty \}$ with the following change of variable: 
\begin{equation}\label{E:transformation variables}
y=-x(x-\epsilon)\frac{w'(x)}{w(x)}
\end{equation}
\end{proposition}

The space of all nonzero solutions ($C_i w_i(x)+ C_j w_j(x)$) of the hypergeometric equation is the manifold $\mathbb{CP}^1 \times \, \mathbb{C}^*$. The next proposition give the expression of a first integral of the Riccati system which takes values in $\mathbb{CP}^1$. Up to a constant (in $\mathbb{C}^*$), this first integral is related to a general solution of the hypergeometric equation.

\begin{proposition}
Let $w_j(X)$ et $w_i(X)$ be two linearly independent solutions of the hypergeometric equation (\ref{E:hypergeo e}). In their shared region of validity we have the following first integral of the Riccati system (\ref{E:ric 1}):
\begin{equation}\label{E:first int new form e}
I^\epsilon_{(i,j)}= \frac{w_i(x)}{w_j(x)} \left( \frac{y  - \rho_i(x,\epsilon)}{y  - \rho_j(x,\epsilon)} \right)
\end{equation}
where
\begin{equation}\label{E:notation rho}
\rho_i(x,\epsilon) =  -x(x-\epsilon )\frac{w_i'(x)}{w_i(x)}.
\end{equation}
\end{proposition}

In order that the limit exists when $\epsilon \in S^+$ goes to zero, we consider the first integral 
\begin{equation}
I^{\epsilon^±}=
\begin{cases}
\kappa^+(\epsilon) I^\epsilon_{(2,3)}  \, \mbox{if } \epsilon \in S^+ \\
\kappa^-(\epsilon) I^\epsilon_{(4,1)} \,  \mbox{if } \epsilon \in S^-
\end{cases}
\end{equation}
where $\kappa^±(\epsilon)$ are defined in (\ref{E:def k}). Now let us see why we can work with a simpler expression than this one to study its ramification.
\begin{proposition}
The quotient $H^{\epsilon^±}=\kappa^±(\epsilon) \frac{w_i(x)}{w_j(x)}$ has the same ramification around $x=0$ and $x=\epsilon$ as 

\begin{equation}
I^{\epsilon^±}= \kappa^±(\epsilon) \frac{w_i(x)}{w_j(x)}\left( \frac{y  - \rho_i(x,\epsilon)}{y  - \rho_j(x,\epsilon)} \right),
\end{equation}
namely we can replace $H^{\epsilon^±}$ by $I^{\epsilon^±}$ in the formulas (\ref{E:H tour e}), (\ref{E:H tour 0}), (\ref{E:K tour e}) and (\ref{E:K tour 0}).

\end{proposition}
\begin{proof}
Let us prove that $H^{\epsilon^+}=\kappa^+(\epsilon) \frac{w_i(x)}{w_j(x)}$ has the same ramification as $I^{\epsilon^+}$ in the case $\epsilon \in S^+$. We start with the ramification around $x=\epsilon$.
We have, with relation (\ref{E:monodromymatrix 2}), 
\begin{equation}\label{E: tour quotient w2}
\begin{array}{lll}
\frac{w'_{2,(\epsilon,-\pi)}(x)}{w_{2,(\epsilon,-\pi)}(x)}&=\frac{\kappa^+(\epsilon) w'_{2,(\epsilon,-\pi)}(x)}{\kappa^+(\epsilon) w_{2,(\epsilon,-\pi)}(x)} \\
&=\frac{ e^{2\pi i (a+b+\frac{1}{\epsilon}-1)}(\kappa^+(\epsilon) w'_{2,(\epsilon,\pi)}(x)-\mu^+(\epsilon) w'_{3,(\epsilon,\pi)}(x))}{ e^{2\pi i (a+b+\frac{1}{\epsilon}-1)}(\kappa^+(\epsilon) w_{2,(\epsilon,\pi)}(x)-\mu^+(\epsilon) w_{3,(\epsilon,\pi)}(x))} \\
&= \frac{1}{ \kappa^+(\epsilon) \frac{w_{2,(\epsilon,\pi)}(x)}{ w_{3,(\epsilon,\pi)}(x)}-\mu^+(\epsilon)}(\kappa^+(\epsilon) \frac{w'_{2,(\epsilon,\pi)}(x)}{ w_{3,(\epsilon,\pi)}(x)}-\mu^+(\epsilon) \frac{w'_{3,(\epsilon,\pi)}(x)}{ w_{3,(\epsilon,\pi)}(x)})\\
&= \frac{1}{ H^{\epsilon^+}_{(\epsilon,\pi)} -\mu^+(\epsilon)} \left( \frac{w'_{2,(\epsilon,\pi)}(x)}{w_{2,(\epsilon,\pi)}(x)} H^{\epsilon^+}_{(\epsilon,\pi)} -\mu^+(\epsilon) \frac{w'_{3,(\epsilon,\pi)}(x)}{w_{3,(\epsilon,\pi)}(x)} \right).
\end{array}
\end{equation}
Using (\ref{E:notation rho}), (\ref{E:H tour e}) and (\ref{E: tour quotient w2}), we have
\begin{equation}
\begin{array}{lll}
I^{\epsilon^+}_{(\epsilon,-\pi)}
&=H^{\epsilon^+}_{(\epsilon,-\pi)}\left( \frac{y  - \rho_{2,(\epsilon,-\pi)}(x,\epsilon)}{y  - \rho_{3,(\epsilon,-\pi)}(x,\epsilon)} \right)\\
&=e^{2 \pi i (a+b-1+\frac{1}{\epsilon})} ( H^{\epsilon^+}_{(\epsilon,\pi)}-\mu^+(\epsilon)) \frac{y +x(x-\epsilon ) \frac{w'_{2,(\epsilon,-\pi)}(x)}{w_{2,(\epsilon,-\pi)}(x)}}{y+ x(x-\epsilon ) \frac{w'_{3,(\epsilon,-\pi)}(x)}{w_{3,(\epsilon,-\pi)}(x)}}\\
&=e^{2 \pi i (a+b-1+\frac{1}{\epsilon})}\frac{ (H^{\epsilon^+}_{(\epsilon,\pi)}-\mu^+(\epsilon ))y + x(x-\epsilon) \left( \frac{w'_{2,(\epsilon,\pi)}(x)}{w_{2,(\epsilon,\pi)}(x)} H^{\epsilon^+}_{(\epsilon,\pi)} -\mu^+(\epsilon) \frac{w'_{3,(\epsilon,\pi)}(x)}{w_{3,(\epsilon,\pi)}(x)} \right) }{y +  x(x-\epsilon ) \frac{w'_{3,(\epsilon,\pi)}(x)}{w_{3,(\epsilon,\pi)}(x)}}\\
&=e^{2 \pi i (a+b-1+\frac{1}{\epsilon})}   \left( H^{\epsilon^+}_{(\epsilon,\pi)}  \frac{y  - \rho_{2,(\epsilon,\pi)}(x,\epsilon)}{y  - \rho_{3,(\epsilon,\pi)}(x,\epsilon)} -\mu^+(\epsilon) \right)  \\
&=e^{2 \pi i (a+b-1+\frac{1}{\epsilon})} \left( I^{\epsilon^+}_{(\epsilon,\pi)}-\mu^+(\epsilon) \right).
\end{array}
\end{equation}
The proofs for $I^{\epsilon^+}_{(0,±\pi)}$, $I^{\epsilon^-}_{(0,±\pi)}$ and $I^{\epsilon^-}_{(\epsilon,±\pi)}$  are similar to this one.
\end{proof}

\subsection{Divergence and unfolding of the saddle-nodes}\label{S:div and unfold}
\begin{figure}[h!]
\begin{center}
{\psfrag{A}{$\small{y=1}$}
\psfrag{B}{$\small{y=0}$}
\psfrag{C}{$\small{x=0}$}
\includegraphics[width=4.5cm]{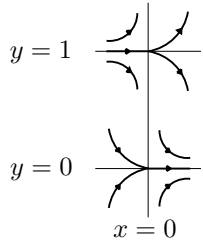}}
    \caption{Phase plane $\epsilon=0$}
    \label{fig:7}
\end{center}
\end{figure}

\begin{figure}[h!]
\begin{center}
{\psfrag{A}{\small{$y=1$}}
\psfrag{B}{\small{$y=y_1$}}
\psfrag{C}{\small{$y=0$}}
\psfrag{D}{\small{$y=0$}}
\psfrag{E}{\small{$x=0$}}
\psfrag{F}{\small{$x=\epsilon$}}
\includegraphics[width=5.5cm]{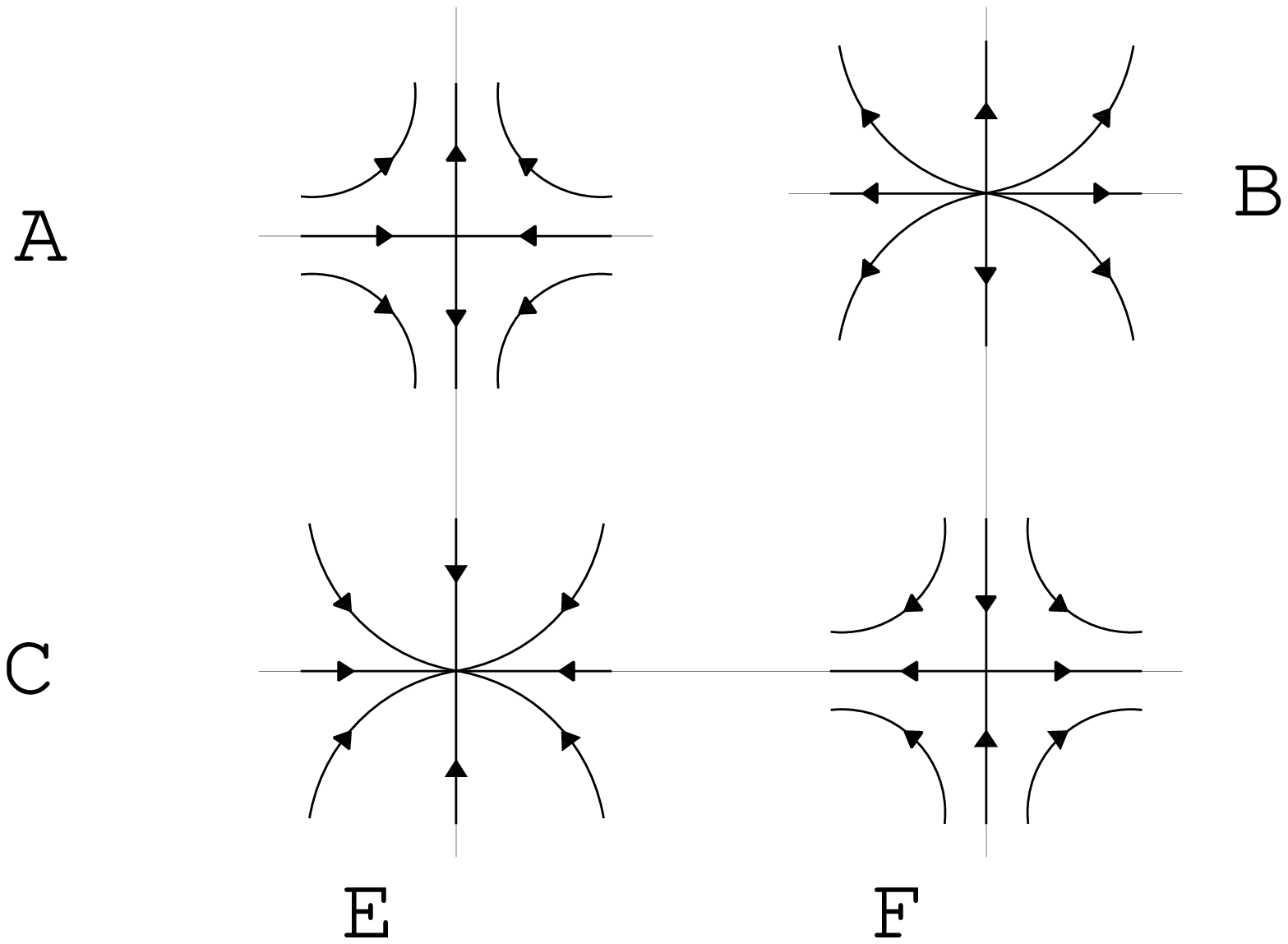}}
    \caption{Phase plane if $\epsilon$ and $\frac{1}{\epsilon}+a+b \in \mathbb{R}$, $\epsilon>0$}
    \label{fig:4}
\end{center}
\end{figure}

Let us consider the Riccati system (\ref{E:ric 1}) with $\epsilon=0$. It has two saddle-nodes located at $(0,0)$ and $(0,1)$ (see Figure \ref{fig:7}). 
In the unfolding (with maybe $a(\epsilon)$ and $b(\epsilon)$), this yields the Riccati system (\ref{E:ric 1}) with the four singular points $(0,0)$, $(\epsilon, 0)$, $(0,1)$ and $(\epsilon, y_1)$ as illustrated in Figures \ref{fig:4} and \ref{fig:5}, with $y_1=1+\epsilon(a+b-1)$. 

The quotient of the eigenvalue in $y$ by the eigenvalue in $x$ of the Jacobian, for each singular point, is given in Table \ref{Ta:prem}.

\begin{table}
	\begin{center}
		\begin{tabular}{|l|l|} \hline
		Singular point  & Quotient of eigenvalues \\ \hline
		$(0,0)$ &  $\frac{1}{\epsilon}$ \\ \hline
		$(\epsilon,0)$ &  $1-\frac{1}{\epsilon}-a-b$\\ \hline
		$(0,1)$ &  $\frac{-1}{\epsilon}$ \\ \hline
		$(\epsilon,y_1)$ & $-1+\frac{1}{\epsilon}+a+b$ \\ \hline
		\end{tabular}

\bigskip

		\caption{Quotient of the eigenvalue in $y$ by the eigenvalue in $x$ of the Jacobian for each singular point}
		\label{Ta:prem}
	\end{center}
\end{table}

\begin{figure}[h!]
\begin{center}
{\psfrag{A}{\small{$y=y_1$}}
\psfrag{B}{\small{$y=1$}}
\psfrag{C}{\small{$y=0$}}
\psfrag{D}{\small{$y=0$}}
\psfrag{E}{\small{$x=\epsilon$}}
\psfrag{F}{\small{$x=0$}}
\includegraphics[width=5.5cm]{FigH}}
    \caption{Phase plane if $\epsilon$ and $\frac{1}{\epsilon}+a+b \in \mathbb{R}$, $\epsilon<0$}
    \label{fig:5}
\end{center}
\end{figure}
 
\begin{remark}\label{R:formal}
By summing the quotient of the eigenvalues at the corresponding saddle and node, we get the formal invariant of the saddle-node at $(0,0)$ (resp. at $(0,1)$), which is $1-a-b$ (resp. $a+b-1$). 
\end{remark}

\begin{figure}[h!]
\begin{center}
{\psfrag{A}{\small{$y=1$}}
\psfrag{B}{\small{$y=y_1$}}
\psfrag{C}{\small{$y=0$}}
\psfrag{D}{\small{$y=0$}}
\psfrag{E}{\small{$x=0$}}
\psfrag{F}{\small{$x=\epsilon$}}
\psfrag{G}{\footnotesize{$y=\rho_2(x,\epsilon)$}}
\psfrag{H}{\footnotesize{$y=\rho_3(x,\epsilon)$}}
\includegraphics[width=5.5cm]{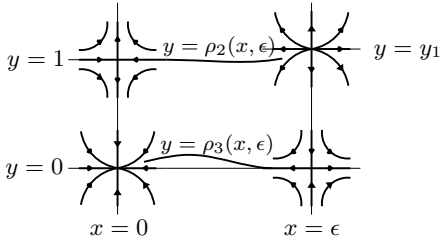}}
    \caption{Invariant manifolds $y=\rho_2(x,\epsilon)$ and $y=\rho_3(x,\epsilon)$, case $\epsilon \in \mathbb{R}^+$}
    \label{fig:8}
\end{center}
\end{figure}

The curves $y-\rho_k(x,\epsilon)=0$ for $k=i,j$ appearing in the first integral (\ref{E:first int new form e}) are solution curves (trajectories) of the Riccati system, more precisely analytic invariant manifolds of two of the singular points when $\epsilon \in S^±$. For example, for $\epsilon \in S^+$, $y=\rho_2(x,\epsilon)$ is the invariant manifold of the singular point $(0,1)$ and $y=\rho_3(x,\epsilon)$ is the invariant manifold of  $(\epsilon,0)$ (see Figure \ref{fig:8}).

Indeed, 
\begin{equation}
\begin{array}{lll}
\rho_2(x,\epsilon) &=- x(x-\epsilon)\frac{w_2'(x)}{w_2(x)} \\
&= 1-\frac{x}{\epsilon}+\{\epsilon(a+b-1)+1\}\frac{x}{\epsilon}+x(1-\frac{x}{\epsilon})\frac{(1-a)(1-b)}{1+\frac{1}{\epsilon}}\frac{\, _2F_1(2-a ,2-b ,2+\frac{1}{\epsilon} ;\frac{x}{\epsilon})}{\, _2F_1(1-a ,1-b ,1+\frac{1}{\epsilon} ;\frac{x}{\epsilon})}
\end{array}
\end{equation}
and $\rho_2(0,\epsilon)=1$. Similarly,
\begin{equation}
\begin{array}{lll}\label{notation rho3}
\rho_3 (x,\epsilon) &= -x(x-\epsilon)\frac{w_3'(x)}{w_3(x)} \\
&=-x(x-\epsilon)\frac{a b}{a+b+\frac{1}{\epsilon}}\frac{\, _2F_1(1+a,1+b,1+a+b+\frac{1}{\epsilon};1-\frac{x}{\epsilon})}{ \, _2F_1(a,b,a+b+\frac{1}{\epsilon};1-\frac{x}{\epsilon})}
\end{array}
\end{equation}
and $\rho_3(\epsilon,\epsilon)=0$.

\begin{figure}[h!]
\begin{center}
\includegraphics[width=4cm]{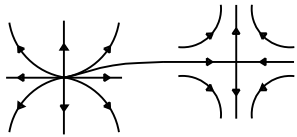}
    \caption{Analytic continuation of an invariant manifold of a saddle when the corresponding analytic center manifold is divergent}
    \label{fig:10}
\end{center}
\end{figure}

The divergence of $g(x)$ corresponds to a nonanalytic center manifold at $(0,0)$ for $\epsilon=0$. When we unfold on $S^+$ (resp. $S^-$), the invariant manifold of $(\epsilon,0)$ (resp. $(0,0)$) is necessarily ramified at $(0,0)$ (resp. $(\epsilon,0)$) for small $\epsilon$ (see Figure \ref{fig:10}). In the particular case when $1-\frac{1}{\epsilon} \in -\mathbb{N}$ (resp. $a+b+\frac{1}{\epsilon}$) with $\epsilon$ small, then $(0,0)$ (resp. $(\epsilon,0)$) is a resonant node. Then necessarily in this case it is non linearisable (the resonant monomial is nonzero) which in practice yields 
logarithmic terms in the first integral. 

\begin{figure}[h!]
\begin{center}
\includegraphics[width=4cm]{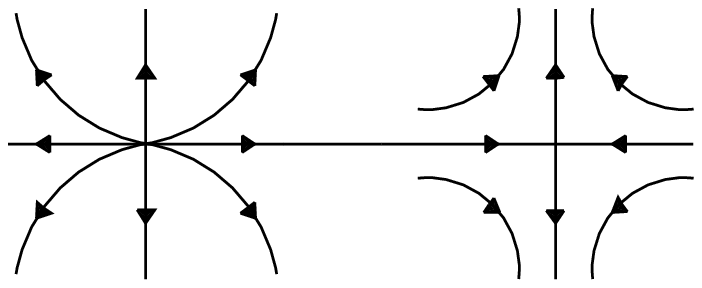}
    \caption{Analytic continuation of an invariant manifold of a saddle when the corresponding analytic center manifold is convergent (this is the case since $a$ and $b$ are fixed)}
    \label{fig:9}
\end{center}
\end{figure}

Besides, if $g(x)$ is convergent, the invariant manifold $y=\rho_3(x)$ (after unfolding in $S^+$, keeping $a$ and $b$ fixed) is not ramified at $(0,0)$ (recall that if $a \in -\mathbb{N}$ or $b \in -\mathbb{N}$, i.e. if $g(x)$ is convergent, then $w_3(x)$ is a polynomial). This correspond to Figure \ref{fig:9}, an exceptional case.

The divergence of $k(x)$ has a similar interpretation with the pair of singular points coming from the unfolding of the saddle-node at $(0,1)$. If $k(x)$ is divergent then, when we unfold in $S^+$ (resp. $S^-$) the invariant manifold of $(0,1)$ (resp. $(\epsilon, y_1)$) is necessarily ramified at $(\epsilon,y_1)$ (resp. $(0,1)$). As before, this implies that $(\epsilon,y_1)$ (resp. $(0,1)$) is nonlinearisable as soon as it is a resonant node. 

The general description of this parametric resurgence phenomenon is described in \cite{cR05}. 

\subsection{Universal unfolding}
As the universal deformation of $x^2$ is $x^2-\epsilon$, let us translate the previous results in the case of this deformation. When studying the universal unfolding of the Riccati system (\ref{E:ric 1}) evaluated at $\epsilon=0$, the singular points to consider would be at $x=-\sqrt{\epsilon}$ and $x=\sqrt{\epsilon}$ (instead of $x=0$ and $x=\epsilon$).

\begin{proposition}
The unfolded Riccati system (with maybe $a(\epsilon)$ and $b(\epsilon)$)
\begin{equation}\label{E:system racine eps}
\begin{cases}
\dot x=x^2-\epsilon\\
\dot y =a(\epsilon) b(\epsilon) (x^2-\epsilon) + (1+(1-a(\epsilon)-b(\epsilon))x)y+y^2
\end{cases}
\end{equation}
is related, with $c=\frac{1}{2\sqrt{\epsilon}}+\frac{a+b+1}{2}$, to the hypergeometric equation with singular points $(-\sqrt{\epsilon},\sqrt{\epsilon},\infty)$
\begin{equation}\label{E:hyp eq racine e}
(x^2-\epsilon) \, w''(x) + \{-1+(a+b+1)x \} \, w'(x)+a b \, w(x) = 0
\end{equation}
with the change of variables
\begin{equation}
y=-(x^2-\epsilon)\frac{w'(x)}{w(x)}.
\end{equation}
\end{proposition}

The product $\lambda^+(\sqrt{\epsilon})\mu^+(\sqrt{\epsilon})$ is an analytic function of $\epsilon$ (and not of $\sqrt{\epsilon}$):

\begin{theorem}
For the family of systems (\ref{E:system racine eps}), in which $a(\epsilon)$ and $b(\epsilon)$ are analytic functions of $\epsilon$, the product $L(\epsilon)=\lambda^+(\sqrt{\epsilon})\mu^+(\sqrt{\epsilon})$ is an analytic function of $\epsilon$. 
\end{theorem}

\begin{proof}
Given $\gamma \in (0,\frac{\pi}{2})$ fixed, we define
\begin{itemize} 
\item $S^+=\{ \epsilon \in \mathbb{C} \, : \, 0<|\epsilon|<r(\gamma), \, \arg(\epsilon) \in (\gamma ,4 \pi -\gamma) \}$.
\end{itemize}

The sector $S^+$ is defined such as $w_2(x)$ and $w_3(x)$ always exist for these values of $\epsilon$. In particular, we ask $-\frac{1}{2\sqrt{\epsilon}}+\frac{3-a+b}{2} \notin -\mathbb{N}$, $-\frac{1}{2\sqrt{\epsilon}}+\frac{a+b+1}{2} \notin -\mathbb{N}$, $-\frac{1}{2\sqrt{\epsilon}}+\frac{a+1-b}{2} \notin -\mathbb{N}$ and $-\frac{1}{2\sqrt{\epsilon}}+\frac{b+1-a}{2} \notin -\mathbb{N}$.

Then, we define
\begin{equation}
H^{\epsilon^+}=\frac{\kappa^+(\sqrt{\epsilon}) w_2(x)}{w_3(x)}
\end{equation}
with 
\begin{equation}
\kappa^+(\sqrt{\epsilon})=(2\sqrt{\epsilon})^{1-a-b}e^{\pi i(\frac{1}{2\sqrt{\epsilon}}+\frac{a+b+1}{2})}
\end{equation}

The functions $\mu^+(\sqrt{\epsilon})$ and $\lambda^+(\sqrt{\epsilon})$ can be defined as before and the calculations give the same relation 
\begin{equation}\label{L(epsilon)}
L(\epsilon)=\lambda^+(\sqrt{\epsilon})\mu^+(\sqrt{\epsilon})=-(1-e^{-2\pi i a(\epsilon)})(1-e^{-2\pi i b(\epsilon)}).
\end{equation}
This product is thus analytic in $\epsilon$ if $a(\epsilon)$ and $b(\epsilon)$ are analytic functions of $\epsilon$.
\end{proof}

These results are used in \cite{cRcC} to characterize the space of modules of a Riccati equation under orbital equivalence. 

\begin{remark}
$L(\epsilon)$ is related to known invariants. Indeed, we have the relation $L(\epsilon)=-4\pi^2 e^{\pi i \alpha(\epsilon)} \gamma(\epsilon) \gamma'(\epsilon)$, where $\alpha(\epsilon)=1-a(\epsilon)-b(\epsilon)$ is the formal invariant of the saddle-node family (\ref{E:system racine eps}), while $\gamma(\epsilon)$ and $\gamma'(\epsilon)$ are the unfolding of the Jurkat-Lutz-Peyerimhoff invariants $\gamma$ and $\gamma'$ (see \cite{jlp}) obtained with the change of variable (\ref{E:chang w-u conf}) in the system associated to the differential equation (\ref{eq hyp confluente originale 2}).
\end{remark}

\section{Directions for further research}
The hypergeometric equation corresponds to a particular Riccati system. The study of this system allowed us to describe how divergence in the limit organizes the system in the unfolding. Similar phenomena are expected to occur in the more general cases where solutions at the confluence are 1-summable or even k-summable. 

\section{Acknowledgements}
We thank the reviewer for useful comments.

\end{document}